\documentclass[a4paper]{article}
\usepackage[dvips]{color,graphicx}
\usepackage{psfrag}
\usepackage{amsmath,amssymb,amsthm,latexsym}

\newtheorem{Thm}{Theorem}
\newtheorem{Prop}{Proposition}
\newtheorem{Cor}{Corollary}

\newtheorem{Definition}{Definition}
\newtheorem*{Example}{Example}
\newtheorem*{Examples}{Examples}
\newtheorem*{Remark}{Remark}
\newtheorem*{Conjecture}{Conjecture}

 \newlength\headseptemp

\newcommand{\Lp}{\Delta}
\newcommand{\Lc}{\widetilde{\Delta}}
\newcommand{\alc}{\widetilde{\al}}

\newcommand{\dd}{{\partial}}

\newcommand{\al}{{\alpha}}
\newcommand{\ka}{{\kappa}}
\newcommand{\ph}{{\varphi}}

\newcommand{\lm}{{\lambda}}

\newcommand{\ra}{{\rightarrow}}

\newcommand{\supp}{{\mathrm{supp}\;}}

\newcommand{\R}{{\mathbb R}}
\newcommand{\D}{{\mathbb D}}

\newcommand{\ab}[1]{\left( #1\right)}

\newcommand{\bs}[1]{\langle #1\rangle}

\newcommand{\ow}[1]{\widetilde{ #1}}

\newcommand{\V}{{\cal V}}
\newcommand{\E}{{\cal E}}
\newcommand{\F}{{\cal F}}
\newcommand{\G}{{\cal G}}
\newcommand{\T}{{\cal T}}
\newcommand{\C}{{\cal C}}
\newcommand{\Cut}{{\rm Cut}}
\newcommand{\vol}{{\rm vol}}

\newcommand{\Hm}[1]{\leavevmode{\marginpar{\tiny%
$\hbox to 0mm{\hspace*{-0.5mm}$\leftarrow$\hss}%
\vcenter{\vrule depth 0.1mm height 0.1mm width \the\marginparwidth}%
\hbox to 0mm{\hss$\rightarrow$\hspace*{-0.5mm}}$\\\relax\raggedright
#1}}}

\date{\today}

\title{Geometric and spectral properties of locally tessellating planar graphs}
\author{Matthias Keller \thanks{e-mail:
    matthias.keller@mathematik.tu-chemnitz.de} \\
  TU Chemnitz\\
  Fakult\"at f\"ur Mathematik\\
  D-09107 Chemnitz, Germany% \vspace*{1ex}
  \and
  Norbert Peyerimhoff \thanks{ e-mail: norbert.peyerimhoff@durham.ac.uk} \\
  Department of Math. Sciences\\
  University of Durham\\
  Durham DH1 2LE, UK}

\begin{document}

\maketitle

%%%%%%%%%%%%%%%%%%%%%%%%%%%%%%%%%%%%%%%%%%%%%%%%%%%%%%%%%%%%
% ABSTRACT
%%%%%%%%%%%%%%%%%%%%%%%%%%%%%%%%%%%%%%%%%%%%%%%%%%%%%%%%%%%%

\begin{abstract} \noindent In this article, we derive bounds for
  values of the global geometry of locally tessellating planar graphs,
  namely, the Cheeger constant and exponential growth, in terms of
  combinatorial curvatures. We also discuss spectral implications for
  the Laplacians.
\end{abstract}

%%%%%%%%%%%%%%%%%%%%%%%%%%%%%%%%%%%%%%%%%%%%%%%%%%%%%%%%%%%%%%%%%%%%%

\section{Introduction}

A locally tessellating planar graph $\G$ is a tiling of the plane with
all faces to be polygons with finitely or infinitely many boundary
edges (see Subsection \ref{loctess} for precise definitions). The edges of
$\G$ are continuous rectifiable curves without self-intersections.
Faces with infinitely many boundary edges are called infinigons and occur, e.g.,
in the case of planar trees. The sets of vertices, edges
and faces of $\G$ are denoted by $\V, \E$ and $\F$. $d(v,w)$ denotes
the combinatorial distance between two vertices $v,w \in \V$, where each
edge is assumed to have combinatorial length one.

Useful {\em local} concepts of the graph $\G$ are combinatorial curvature
notions. The finest curvature notion is defined on the corners of
$\G$. A corner is a pair $(v,f) \in \V \times \F$, where $v$ is a
vertex of the face $f$. The set of all corners is denoted by $\C$. The
{\em corner curvature} $\kappa_C$ is then defined as
$$ \kappa_C(v,f) = \frac{1}{|v|} + \frac{1}{|f|} - \frac{1}{2}, $$
where $|v|$ and $|f|$ denote the degree of the vertex $v$ and the face
$f$. If $f$ is an infinigon, we set $|f| = \infty$ and $1/|f| = 0$.
The {\em curvature at a vertex} $v \in \V$ is given by the sum
$$ \kappa(v) = \sum_{(v,f) \in \C}\kappa_C(v,f) = 1 - \frac{|v|}{2} +
\sum_{f: v \in f} \frac{1}{|f|}. $$ For a finite set $W \subset \V$
we define $\ka(W) = \sum_{v \in W} \ka(v)$. These combinatorial
curvature definitions arise naturally from considerations of the
Euler characteristic and tessellations of closed surfaces, and they
allow to prove a combinatorial Gau{\ss}-Bonnet formula (see
\cite[Thm 1.4]{BP1}). Similar combinatorial curvature notions have
been introduced by many other authors, e.g., \cite{Gro,Hi,St,Woe}.

The aim of this paper is to establish connections between local
curvature conditions and characteristic values of the {\em global geometry}
of the graph $\G$, in particular the exponential growth and Cheeger
constants. For a finite subset $W \subset \V$, let $\vol(W) = \sum_{v
  \in W} |v|$. The exponential growth is defined as follows (note that
the value $\mu(\G)$ does not depend on the choice of center $v \in
\V$):

\begin{Definition} The {\em exponential growth} $\mu(\G)$ is given by
$$ \mu(\G) = \limsup_{n \to \infty} \frac{\log \vol(B_n(v))}{n}, $$
where $B_n(v) = \{ w \in \V \mid d(v,w) \le n \}$ denotes the
(combinatorial) ball of radius $n$ about $v$.
\end{Definition}

We also consider the following two types of Cheeger constants.

\begin{Definition} \label{d:cheeg}
  Let
  \begin{equation*}
    \al(\G)=\inf_{\scriptsize\begin{array}{c}
        W\subseteq \V, \\
        |W|<\infty
      \end{array}}
    \frac{|\dd_E W|}{|W|}\quad\mathrm{and}\quad
    \alc(\G)=\inf_{\scriptsize\begin{array}{c}
        W\subseteq \V, \\
        |W|<\infty
      \end{array}}
    \frac{|\dd_E W|}{\vol(W)}
  \end{equation*}
  where $\dd_E W$ is the set of all edges $e \in \E$ connecting a
  vertex in $W$ with a vertex in $\V \backslash W$. $\al(\G)$ is
  called the \emph{physical Cheeger constant} and $\alc(\G)$ the
  \emph{combinatorial Cheeger constant} of the graph $\G$.
\end{Definition}

The attributes {\em physical} and {\em combinatorial} in the previous
definition are motivated by the fact that these Cheeger constants are
closely linked to two types of Laplacians: The {\em physical
  Laplacian} is used frequently in the community of Mathematical
Physicists and is defined as follows:
\begin{equation} \label{laplp}
(\Lp\ph)(v)=|v|\ph(v) - \sum_{w\sim v}\ph(w).
\end{equation}
Note that $\Lp$ is an unbounded operator if there is no bound on the
vertex degree of $\G$. The {\em combinatorial Laplacian} $\Lc$ is a
bounded operator and appears in the context of spectral geometry
(see, e.g., \cite{DKa,DKe,Woe2}):
\begin{equation} \label{laplc}
(\Lc\ph)(v)=\ph(v) - \frac{1}{|v|}\sum_{u\sim v}\ph(u).
\end{equation}
Both operators are defined in and are self-adjoint with respect to
different $l^2$-spaces (see Subsection \ref{lap}). In the case of
fixed vertex degree, both operators are multiples of each
other.\medskip

Our main geometric results are given in Subsection \ref{cheegexp},
where we
\begin{itemize}
\item provide lower bounds for both Cheeger constants in terms of
  combinatorial curvatures (see Theorem \ref{t:cheeg} below),
\item provide upper bounds for the exponential growth in terms of an upper
  vertex bound (see Theorem \ref{t:expgr} below).
\end{itemize}
Even though Theorem \ref{t:expgr}(b) is formulated in terms of bounds
on vertex and face degrees, it can also be considered as an estimate
in terms of combinatorial curvature, as is explained in the remark
following the theorem. In fact, the proof is based on the
corresponding curvature version.

\medskip

Now we discuss connections to the spectrum. The Cheeger constant and
the exponential growth were first introduced in the context of
Riemannian manifolds and were useful invariants to estimate the bottom
of the (essential) spectrum of the Laplacian (see \cite{Che} and
\cite{Br}). An analogous inequality between the Cheeger constant and
the bottom of the spectrum in the discrete case of graphs was first
proved by \cite{Do} and \cite{Al}. This inequality is also useful in
the study of expander graphs. \cite{Al} noted also the connection
between this inequality and the Max Flow-Min Cut Theorem (see also
\cite{Chu} and \cite{Gri}). For other connections between isoperimetric
inequalities and lower bounds of eigenvalues in both continuous and
discrete settings see, e.g., \cite{CGY}.

The best results about the relations between the combinatorial Cheeger
constant, the exponential growth, and the bottom $\ow\lm_0(\G)$ and
$\ow\lm_0^{ess}(\G)$ of the (essential) spectrum of the {\em combinatorial
Laplacians} $\Lc$ are due to K. Fujiwara (see \cite{Fu1} and
\cite{Fu2}):
\begin{equation} \label{fuji}
1-\sqrt{1-\widetilde\alpha^2(\G)} \le \ow\lambda_0(\G) \le \ow\lambda_0^{ess}(\G)
\le 1 - \frac{2 e^{\mu(\G)/2}}{1+e^{\mu(\G)}}.
\end{equation}
These estimates are sharp in the case of regular trees. Using these
estimates and Theorems \ref{t:cheeg} and \ref{t:expgr}, we obtain
\begin{itemize}
\item lower and upper estimates on the bottom of the (essential)
spectrum of the combinatorial Laplacian in terms of combinatorial
curvature (see Corollaries \ref{c:McKean} and \ref{c:essenspecest}).
\end{itemize}
Since there are estimates to compare the bottom of the (essential)
spectrum of the combinatorial Laplacian with the physical Laplacian
(see for instance \cite{Ke}) these results can be also formulated
for the physical Laplacian.

A lower estimate for the bottom of the essential spectrum of the
combinatorial Laplacian via the combinatorial Cheeger constant at
infinity can be found in \cite[Cor. 3]{Fu2}. This yields a discrete
analogue for the combinatorial Laplacian of the result in \cite{DL}
about the {\em emptiness of the essential spectrum} for complete simply
connected manifolds with curvature converging to minus infinity.
Corresponding results about the emptiness of the essential spectrum
for the physical Laplacian can be found in \cite{Ke,Woj}.

Finally, let us discuss two other interesting types of eigenfunctions,
namely, {\em strictly positive eigenfunctions} and {\em finitely supported
eigenfunctions}, and illustrate all concepts in two examples.

For the discrete case of a graph, it was shown in
\cite[Prop. 1.5]{DKa} that the equation $\Lc f = \lambda f$ has a
{\em positive solution} if and only if $\lambda \le \ow\lambda_0(\G)$. This
characterisation of the bottom of the spectrum was well known before in the
context of Riemannian manifolds (see, e.g., \cite{Sull} and the
references therein). In the reverse direction, this characterisation
might be used in concrete cases to determine the bottom of the spectrum
of an infinite graph.

On the other hand, {\em finitely supported solutions} of the equation
$\Lc f = \lambda f$ are obviously $l^2$-eigenfunctions and, therefore,
they can only exist for eigenvalues $\lambda \ge
\ow\lambda_0(\G)$. Existence of finitely supported eigenfunctions in
Penrose tilings was first observed in \cite{KS}. Their
existence is a purely discrete phenomenon, since in the case of a
non-compact, connected Riemannian manifold the eigenvalue equation
$\Delta f = \lambda f$ cannot have compactly supported eigenfunctions
(a fact which is known as the {\em unique continuation principle}; see
\cite{Ar}). These finitely supported eigenfunctions coincide with the
discontinuities of the integrated density of states (or spectral
density function). See, e.g., the articles \cite{KLS,LV} and the references
therein for more details about this connection.

\begin{Examples}
(a) We consider the periodic tessellation $\G = (\V,\E,\F)$ in
  Figure \ref{example1}. We assume that all edges are straight Euclidean
  segments of length one.

  \begin{figure}[h]
      \psfrag{v}{$v$}
      \begin{center}
      \includegraphics[height=6cm]{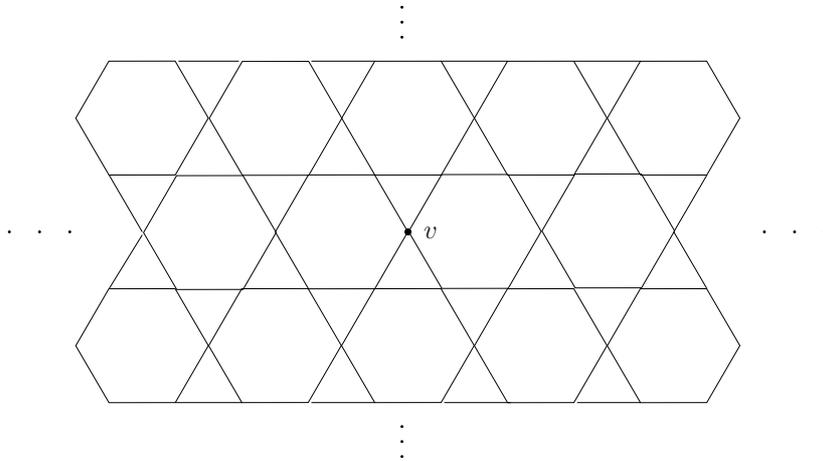}
    \end{center}
    \caption{Plane tessellation with regular triangles and hexagons}
    \label{example1}
  \end{figure}

  We first show that $\mu(\G) =0$: Choose a fixed radius
  $0<r<1/2$. Then all Euclidean balls of radius $r$ centered at all
  vertices in $\V$ are pairwise disjoint. On the other hand, the
  vertices in the combinatorial ball $B_n(v)$ are contained in the
  Euclidean ball of radius $n$, centered at $v$. Both facts together
  imply that combinatorial balls grow only polynomially and the
  exponential growth is zero. As a consequence, this graph cannot
  contain a binary tree as a subgraph. Moreover, using \eqref{fuji}, we conclude that
  $$ \ow\lambda_0(\G) = \ow\lambda_0^{ess}(\G) =0\quad \mathrm{and} \quad \alc(\G) = \al(\G)
  = 0. $$
  Finally, $\G$ does admit finitely supported eigenfunctions, namely,
  choose $p \in \R^2$ to be the center of a hexagon and define $f(p +
  e^{2\pi i/6}) = (-1)^i$ (i.e., choose alternating values
  $1,-1,1,-1,1,-1$ clockwise around the vertices of the hexagon) and
  $f(v) = 0$ for all other vertices. Then we have $\Lc f =\frac{3}{2}
  f$.\medskip

(b) Let $\T_p$ denote the $p$-regular tree. In this case,
  spectrum and essential spectrum of the combinatorial Laplacian
  coincide and are given by the interval (see, e.g., \cite[App. 3]{Sun})
  $$\left[1-\frac{2\sqrt{p-1}}{p}, 1+\frac{2\sqrt{p-1}}{p}\right].$$
  Consequently, $\Lc f = \lambda f$ admits a positive solution if and
  only if $\lambda \le 1-2\sqrt{p-1}/p$. Moreover, we have $\alc(\T_p)
  = \frac{p-2}{p}$, $\al(\T_p) = p-2$ and $\mu(\T_p) = \log
  (p-1)$. Note that a regular tree doesn't admit
  $l^2$-eigenfunctions. For otherwise, we could choose a vertex $v$ at
  which our eigenfunction doesn't vanish and take its radialisation
  with respect to this vertex. This radialisation would be again a
  non-vanishing $l^2$-eigenfunction with the same eigenvalue and,
  since its values would only depend on the distance to $v$, there
  would be an easy recursion formula for its values. The precise form
  of the recursion formula would then contradict to the requirement
  that the function lies in $l^2$.
\end{Examples}
\textbf{Acknowledgements.} Matthias Keller would like to
thank Daniel Lenz who encouraged him to study the connection between
curvature and spectral theory. Matthias Keller was supported during
this work by the German Business Foundation (sdw).

%%%%%%%%%%%%%%%%%%%%%%%%%%%%%%%%%%%%%%%%%%%%%%%%%%%%%%%%%%%%%%%%%%%%%

\section{Basic notions and main results}

In the first two subsections, we provide the notions which
haven't yet been introduced in full detail in the Introduction. In
Subsections \ref{cheegexp} and \ref{specappl}, we state our main results.

%%%%%%%%%%%%%%%%%%%%%%%%%%%%%%%%%%%%%%%%%%%%%%%%%%%%%%%%%%%%%%%%%%%%%%%%%%%

\subsection{Locally tessellating planar graphs}
\label{loctess}

Let $\G=(\V,\E)$ be a planar graph (with $\V$ and $\E$ the set of
vertices and edges) embedded in $\R^2$. The faces $f$ of $\G$ are the
closures of the connected components in $\R^2\setminus \bigcup_{e\in
  E} e$. The set of faces is denoted by $\F$.

We further assume that $\G$ has no loops, no multiple edges and no
vertices of degree one (terminal vertices). We write $e = vw$, if the
edge $e$ connects the vertices $v,w$. Moreover, we assume that every
vertex has finite degree and that every bounded open set in $\R^2$
meets only finitely many faces of $\G$. We call a planar graph with
these properties {\em simple}. The {\em boundary of a face $f$} is the
subgraph $\partial f = (\V \cap f,\E \cap f)$. We call a sequence of
edges $e_1,\dots,e_n$ a {\em walk of length $n$} if there is a
corresponding sequence of vertices $v_1,\dots,v_{n+1}$ such that $e_i
= v_iv_{i+1}$. A walk is called a {\em path} if there is no repetition
in the corresponding sequence of vertices $v_1,\dots,v_n$.

A simple planar graph $\G$ is called a {\em locally
tessellating planar graph} if the following additional conditions are
satisfied:
\begin{itemize}
\item [i.)] Any edge is contained in precisely two different faces.
\item [ii.)] Any two faces are either disjoint or have precisely a
  vertex or a path of edges in common. In the case that the length of
  the path is greater then one, then both faces are unbounded.
\item [iii.)] Any face is homeomorphic to the closure of an open disc
  $\D \subset \R^2$, to $\R^2\setminus \D$ or to the upper half plane
  $\R\times\R_+\subset \R^2$ and its boundary is a path.
\end{itemize}
Note that these properties force the graph $\G$ to be connected.
Examples are tessellations $\R^2$ introduced in \cite{BP1,BP2}, trees
in $\R^2$, and particular finite tessellations on the sphere mapped to
$\R^2$ via stereographic projection.

When we consider the vertex degree as a function on $\V$ we write
$\deg(v)=|v|$ for $v\in \V$.  Moreover we define the degree $|f|$ of
a face $f\in F$ to be the length of the shortest closed walk in the
subgraph $\partial f$ meeting all its vertices. If there is no such
finite walk we set $|f|=\infty$. $v \sim w$ means that $d(v,w)=1$,
i.e., $v$ and $w$ are neighbors. A (finite or infinite) path with
associated vertex sequence $\dots v_i v_{i+1} v_{i+2} \dots$ is
called a {\em geodesic}, if we have $d(v_i,v_j) = |i-j|$ for all
pairs of vertices in the path.

%%%%%%%%%%%%%%%%%%%%%%%%%%%%%%%%%%%%%%%%%%%%%%%%%%%%%%%%%%%%%%%%%%%%%%%%%%%

\subsection{Laplacians}
\label{lap}

Let $\G = (\V,\E,\F)$ be a locally tessellating planar graph.
The operators $\Lp$ and $\Lc$ were already introduced in \eqref{laplp}
and \eqref{laplc}. They are symmetric operators and initially defined
on the space
$$ c_c(\V) := \{\ph: \V \ra \R \mid\, |\supp\ph| <\infty \} $$
of functions with finite support. However, they have unique self-adjoint
extensions on different $l^2$-spaces: Let $g: \V \to (0,\infty)$ be a
weight function on the vertices of the graph $\G$ and
$$ l^2(\V,g) := \{\ph:\V\ra\R\mid \bs{\ph,\ph}_g:=\sum_{v\in V}
g(v)|\ph(v)|^2<\infty\}. $$
For $g=1$ we simply write $\l^2(\V)$.

Then the combinatorial Laplacian can be extended to a bounded
self-adjoint operator on all of $l^2(\V,\deg)$. The physical Laplacian
has also a unique self-adjoint extension in the space $l^2(\V)$ (see
\cite{We} or \cite{Woj}). Note, however, that the adjacency operator
need not be essentially self adjoint (see \cite[Section 3]{MW} and the
references therein). We denote the self-adjoint extensions of both Laplacians,
again, by $\Lc$ and $\Lp$.

Furthermore, we define the restriction of the combinatorial Laplacian on the
complement of a finite set $K$ of vertices. Let $P_K:l^2(\V,\deg)\ra
l^2(\V\setminus K,\deg)$ be the canonical projection and
$i_K:l^2(\V\setminus K,\deg)\ra l^2(\V,\deg)$ be its dual operator,
which is the continuation by $0$ on $K$. We write $\Lc_K=P_K \Lc
i_K$. Of particular importance is the bottom of the spectrum
$\ow\lambda_0(\G)$ and of the essential spectrum
$\ow\lambda_0^{ess}(\G)$. $\ow\lambda_0(\G)$ can be characterised as the
infimum of the Rayleight-Ritz quotient over all non-zero functions $f
\in l^2(\V,\deg)$, i.e.,
$$
\ow\lambda_0(\G) = \inf \left\{ \frac{\langle \Lc f,f\rangle_{\deg}}{\langle
  f,f\rangle_{\deg}}: f \neq 0, f \in l^2(\V,\deg) \right\}.
$$
Similarly, $\ow\lambda_0^{ess}(\G)$ can be obtained via
\begin{equation} \label{specess}
\ow\lambda_0^{ess}(\G) = \lim_{n \to \infty} \inf \left\{ \frac{\langle
    \Lc_{B_n} f,f\rangle_{\deg}}{\langle f,f\rangle_{\deg}}: f \neq 0, f
  \in l^2(\V \backslash B_n,\deg) \right\},
\end{equation}
where $B_n$ are balls of radius $n$ around any fixed vertex $v \in
\V$. A proof of \eqref{specess} can be found in \cite{Ke}. Obviously, we have
$\ow\lambda_0(\G) \le \ow\lambda_0^{ess}(\G)$. Equality holds in the
following case:

\begin{Prop} \label{specessspeceq}
  Assume that there is a subgroup $\Gamma$ of the automorphism group
  of $\G$ with $\sup_{\gamma \in \Gamma} d(v,\gamma v) = \infty$ for some
  vertex $v \in \V$. Then we have
  $$ \ow\lambda_0(\G) = \ow\lambda_0^{ess}(\G). $$
\end{Prop}

\begin{proof}
For the bottom of the spectrum not to lie in the essential spectrum
would mean that it is an isolated eigenvalue of finite multiplicity.
But this cannot be the case (see Fact 1 in \cite[p. 259]{Sun}).
\end{proof}

Analogous statements hold for the bottom of the (essential) spectrum of
the physical Laplacian.

%%%%%%%%%%%%%%%%%%%%%%%%%%%%%%%%%%%%%%%%%%%%%%%%%%%%%%%%%%%%%%%%%%%%%%%%%%%

\subsection{Cheeger constant and exponential growth estimates}
\label{cheegexp}

The physical and combinatorial Cheeger constants were introduced in
Definition \ref{d:cheeg}. It is easy to see that they are linked to the
physical and combinatorial Laplacians via the equations:
\begin{equation*}
  \al(\G)=\inf_{\scriptsize\begin{array}{c}
      W\subseteq \V, \\
      |W|<\infty
    \end{array}}
  \frac{\bs{\Lp\chi_W,\chi_W}}{\bs{\chi_W,\chi_W}}\quad\mathrm{and}\quad
  \alc(\G)=\inf_{\scriptsize\begin{array}{c}
      W\subseteq \V, \\
      |W|<\infty
    \end{array}}
  \frac{{\bs{\Lc\chi_W,\chi_W}}_{\deg}}{{\bs{\chi_W,\chi_W}}_{\deg}},
\end{equation*}
where $\chi_W$ denotes the characteristic function of the set
$W\subseteq \V$. Note, in particular, that the combinatorial Cheeger
constant is always bounded from above by $\alc(\G) \le 1$.

Next, we state the Cheeger constant estimates:

\begin{Thm} \label{t:cheeg}
  Let $\G =(\V,\E,\F)$ be a locally tessellating planar graph and $3
  \le q \le \infty\}$ such that $|f|\leq q$ for all faces $f\in F$.
  \begin{itemize}
  \item[(a)] For some $a > 0$, let $\ka(v) \le -a$ for all $v \in \V$.
  Then we have
    $$\al(\G) \geq \frac{2q}{q-2}a.$$
  \item[(b)] For some $c > 0$, let
    $\frac{1}{|v|}\ka(v) \le -c$ for all $v \in
    \V$. Then we have
    $$\alc(\G) \geq \frac{2q}{q-2} c.$$
  \end{itemize}
  Moreover, the above estimates are sharp in the case of regular
  trees. (Note that in the case $q=\infty$ we set $\frac{2q}{q-2} = 2$.)
\end{Thm}

\begin{Remark}
  The combinatorial Cheeger constant of all non-positively curved {\em
    regular} plane tessellation $\G_{p,q}$ (with all vertices
  satisfying $|v| = p$ and faces satisfying $|f| = q$) was explicitly
  calculated in \cite{HJL} and \cite{HiShi} as
  $$ \alc(\G_{p,q}) = \frac{p-2}{p} \sqrt{1 - \frac{4}{(p-2)(q-2)}}. $$
  Our estimate gives in this case
  $$ \alc(\G_{p,q}) \ge \frac{(p-2)(q-2)-4}{p(q-2)}. $$
\end{Remark}

Before considering the exponential growth of a locally tessellating
planar graph $\G=(\V,\E,\F)$, let us first introduce the {\em cut
  locus} $\Cut(v)$ of a vertex $v \in \V$. $\Cut(v)$ denotes the set
of all vertices $w$, at which $d_v := d(v,\cdot)$ attains a local
maximum, i.e., we have $w \in \Cut(v)$ if $d_v(w') \le d_v(w)$ for all
$w' \sim w$. $\G$ is {\em without cut locus} if $\Cut(v) = \emptyset$
for all $v \in \V$. Obviously, the cut locus of a finite graph is
never empty. It was proved in \cite[Thm. 1]{BP2} that plane
tessellations with everywhere non-positive corner curvature are graphs
without cut locus. Moreover, let $\T_p$ denote the regular tree with
$|v| = p$ for all vertices.

\begin{Thm} \label{t:expgr}
  Let $\G =(\V,\E,\F)$ be a locally tessellating planar graph without
  cut locus.
  \begin{itemize}
  \item[(a)] If there exists $p \ge 3$ such that
  \begin{equation} \label{vertexest}
  | v | \le p \quad \forall\, v \in \V,
  \end{equation}
  then we have
  $$ \mu(\G) \le \mu(\T_p) = \log(p-1). $$
  \item[(b)] If there exist $p \ge 3$ such that \eqref{vertexest} is
  satisfied and $q \in \{3,4,6\}$ such that
  $$ | f | = q \quad \forall\, f \in \F, $$
  (i.e., $\G$ is face-regular) then we have
  $$ \mu(\G) \le \mu(\G_{p,q}) = \log \left( \frac{p}{2} -
    \frac{2}{q-2} + \sqrt{ \left( \frac{p}{2} - \frac{2}{q-2}
      \right)^2 - 1 } \right). $$
  \end{itemize}
\end{Thm}

\begin{Remark}
  For the reader's convenience, Theorem \ref{t:expgr}(b) was stated in
  ``more familiar'' terms of vertex and face degrees. However, the
  statement has an equivalent reformulation in terms of curvature: Let
  $\G$ be a locally tessellating planar graph without cut locus
  satisfying $|f| = q$ for all faces and $q \in \{3,4,6\}$. For some
  $b \ge 0$, let $-b \le \ka(v)$ for all $v \in \V$. Then we have
  $$ \mu(\G) \le \log ( \tau + \sqrt{\tau^2 -1} ), $$
  where $\tau = 1 + \frac{q}{q-2} b \ge 1$. The inequality is sharp
  (with the optimal choice of $b$) in the case of regular graphs
  $\G_{p,q}$. In fact, the proof will be given for this equivalent
  reformulation. (Note that the constants $p$ and $b$ in the two
  formulations are related by $b = \frac{q-2}{q} p - 1$.)
\end{Remark}

Since the regular plane tessellations $\G_{p,q}$ can be considered as
combinatorial analogues of {\em constant curvature space forms} in
Riemannian geometry, it is natural to conjecture the following
discrete version of a {\em Bishop volume comparison result} (see,
e.g., \cite[Theorem 3.101]{GaHuLa} for the case of a Riemannian
manifold).

\begin{Conjecture}
  Let $p,q \ge 3$ with $1/p+1/q \le 1/2$ be given. Then we have
  \begin{equation} \label{bishop} \mu(\G) \le \mu(\G_{p,q}),
  \end{equation}
  for all locally tessellating planar graphs $\G = (\V,\E,\F)$ without
  cut locus satisfying $|v| \le p$, $|f| \le q$.
\end{Conjecture}

Theorem \ref{t:expgr} confirms this conjecture for the cases $q = 3$
and $q = \infty$. However, it seems difficult to prove this seemingly
obvious estimate \eqref{bishop} for general face degree bounds $q \ge
3$. Assuming the above conjecture to be true, the comparison of the
exponential growth of a locally tessellating planar graph with upper
vertex degree bound $p$ and of the regular tree $\T_p$, as given in
Theorem \ref{t:expgr}(a), is quite good if all faces of $\G$ satisfy
$|f| \ge 6$. For example, we have in the case $(p,q)=(5,6)$:
$$ 1.307\dots = \log(2+\sqrt{3}) = \mu(\G_{5,6}) \le \mu(\T_5) = \log 4 =
1.381\dots. $$

An direct consequence of \cite[Corollary 5.2]{BP1} is the following
lower bound for the exponential growth:

\begin{Thm} \label{t:bp}
  Let $\G = (\V,\E,\F)$ be a locally tessellating planar graph without
  cut locus and $a > 0$ such that $\ka(v) \le -a$ for all vertices
  $v \in \V$. Assume there is $3 \le q \le \infty$ such
  that we have $|f| \le q$ for all faces $f \in \F$. Then we have
  $$ \mu(G) \ge \log\left( 1 + \frac{2q}{q-1} a \right). $$
  Moreover, this estimate is sharp in the case of regular trees. (In
  the case $q = \infty$, we set $\frac{2q}{q-1} = 2$.)
\end{Thm}

We like to finish this subsection by a few additional useful facts: Let
$$ S_n(v) =  \{ w \in \V \mid d(v,w) = n \} $$
be the (combinatorial) sphere of radius $n$ about $v \in \V$. If there
is a uniform upper bound on the vertex degree and if $s_n := | S_n(v)
|$ is a non-decreasing sequence, one easily checks that
\begin{equation} \label{muGalt} \mu(\G) = \limsup_{n \to \infty}
  \frac{\log s_n}{n}.
\end{equation}
Yet another Cheeger constant $h(\G)$ was considered in \cite{BS}:
$$ h(\G) =  \inf_{\scriptsize\begin{array}{c}
      W\subseteq \V, \\
      |W|<\infty
      \end{array}}
    \frac{|\dd_V W|}{|W|},
$$
where $\dd_V W$ is the set of all vertices $v \in \V \backslash W$
which are end points of an edge in $\dd_E W$. In the case that
$\mu(\G)$ is presented by \eqref{muGalt}, this Cheeger constant is
related to the exponential growth by
$$ e^{\mu(\G)} \ge 1+h(\G), $$
with equality in the case of regular trees.

%%%%%%%%%%%%%%%%%%%%%%%%%%%%%%%%%%%%%%%%%%%%%%%%%%%%%%%%%%%%%%%%%%%%%%%%%%%

\subsection{Spectral applications}
\label{specappl}

An immediate consequence of Fujiwara's lower estimate \eqref{fuji}
and Theorem \ref{t:cheeg} is the following {\em combinatorial analogue
of McKean's Theorem} (see \cite{McK} for the case of a Riemannian
manifold):

\begin{Cor}[Combinatorial version of McKean's
  Theorem] \label{c:McKean} Let $\G = (\V,\E,\F)$ be a locally
  tessellating planar graph and $3 \le q \le \infty$ such that
  $|f|\leq q$ for all faces $f\in F$. For some $c > 0$, let
  $\frac{1}{|v|} \ka(v) \le -c$ for all $v \in \V$. Then we have
  $$ 1-\sqrt{1-\left( \frac{2q}{q-2} c\right)^2} \le \ow\lambda_0(\G). $$
  This estimate is sharp in the case of regular trees.
\end{Cor}

Combining Theorem \ref{t:expgr}(a), the curvature version of Theorem
\ref{t:expgr}(b) (see the remark of the theorem) and Fujiwara's upper
estimate \eqref{fuji}, we obtain:

\begin{Cor} \label{c:essenspecest}
  Let $\G =(\V,\E,\F)$ be a locally tessellating planar graph without
  cut locus.
  \begin{itemize}
  \item[(a)] If there exists $p \ge 3$ such that
    \begin{equation} \label{vertexest2} | v | \le p \quad \forall\, v
      \in \V,
    \end{equation}
    then we have
    $$
    \ow\lambda_0^{ess}(\G) \le \ow\lambda_0^{ess}(\T_p) = 1 -
    \frac{2\sqrt{p-1}}{p}.
    $$
  \item[(b)] If there exist $q \in \{3,4,6\}$ with $|f| = q$ for all
  $f \in \F$, and $b > 0$ with $-b \le \ka(v)$ for all $v \in \V$, then
  we have
    $$
    \ow\lambda_0^{ess}(\G) \le 1 - \frac{2 \sqrt{\tau + \sqrt{ \tau^2 - 1
        }} }{1 + \tau + \sqrt{ \tau^2 - 1 } },
    $$
  where $\tau = 1 + \frac{q}{q-2}b$.
  \end{itemize}
\end{Cor}

Next we indicate implications of the above results for the spectrum
of the {\em physical Laplacian}.  Let $\lm_0(\G)$ and
$\lm_)^{ess}(\G)$ denote the bottom of the (essential) spectrum of
the physical Laplacian $\Delta$ and, for $n \ge 0$, let
$$m_n=\inf_{w\in V\setminus B_{n-1}(v)} |w| \quad \mathrm{and}
\quad M_n=\sup_{w\in V\setminus B_{n-1}(v)} |w|,$$ where $v \in \V$
is an arbitrary vertex and $B_{-1}(v) = \emptyset$. Moreover let
$m_\infty=\lim_{n \to \infty} m_n$ and $M_\infty=\lim_{n \to \infty}
M_n$. Then we have, by \cite{Do}
\begin{equation} \label{dodziuk}
\lm_0(\G) \ge \frac{\al(\G)^2}{2M} \quad \mathrm{and} \quad
\lm_0^{ess}(\G) \ge \frac{\al_\infty(\G)^2}{2M_\infty},
\end{equation}
where $\al_\infty(\G)$ denotes the physical Cheeger constant at
infinity, defined in \cite{Ke}. In general we can also estimate, as
demonstrated in \cite{Ke},
$$
m_0 \ow\lm_0(\G) \leq \lm_0(\G) \leq M_0 \ow\lm_0(\G) \quad
\mathrm{and} \quad m_\infty \ow\lm_0^{ess}(\G) \leq \lm_0^{ess}(\G)
\leq M_\infty \ow\lm_0^{ess}(\G).
$$
Via this inequalities we can estimate the bottom of the (essential)
spectrum of the physical Laplacian $\Lp$ by the estimates of
Corollary \ref{c:McKean} and \ref{c:essenspecest} for the
combinatorial Laplacian.

Before we look at an explicit example, let us mention the following
result about the absence of finitely supported eigenfunctions in the
case of non-positive corner curvature:

\begin{Thm}[see {\cite[Theorem 4]{KLPS}}] \label{t:KLPS} Let $\G =
  (\V,\E,\F)$ be a plane tessellation (in the restricted sense of
  \cite{BP2}) with non-positive corner curvature in all corners. Then
  the combinatorial Laplacian does not admit finitely supported
  eigenfunctions.
\end{Thm}

Note that Theorem \ref{t:KLPS} becomes wrong if we replace
``non-positive corner curvature'' by the weaker assumption
``non-positive vertex curvature'', since Example (a) of the
Introduction is a graph with vanishing vertex curvature which admits
finitely supported eigenfunctions.

\medskip

Let us, finally, apply the above results in an example.

\begin{Example}
  We consider the regular tessellation $\G_{6,6}$. Using our geometric
  results in this article, we obtain
  $$ \alc(\G_{6,6}) \ge \frac{1}{2} \quad \textrm{and} \quad \mu(\G_{6,6}) =
  \log \frac{1+\sqrt{21}}{2} \approx 1.5668. $$
  Proposition \ref{specessspeceq} tells us that $\ow\lambda_0(\G_{6,6}) =
  \ow\lambda_0^{ess}(\G_{6,6})$, and with our results in this Subsection
  we can conclude that
  \begin{eqnarray*}
    \ow\lambda_0(\G_{6,6}) = \ow\lambda_0^{ess}(\G_{6,6}) &\in& \left[
      1-\frac{\sqrt{3}}{2}, 1- 2 \frac{\sqrt{3}+\sqrt{7}}{7+\sqrt{21}}
    \right] \\
    &\approx& [ 0.1340, 0.2441 ].
  \end{eqnarray*}
  Using the explicit formula for the Cheeger constant in \cite{HJL} in
  this particular case, we obtain $\alc(\G_{6,6}) = \frac{1}{\sqrt{3}}
  \approx 0.5774$ and we can shrink this interval to
  \begin{eqnarray*}
    \ow\lambda_0(\G_{6,6}) = \ow\lambda_0^{ess}(\G_{6,6}) &\in& \left[
      1-\sqrt{\frac{2}{3}}, 1- 2 \frac{\sqrt{3}+\sqrt{7}}{7+\sqrt{21}}
    \right] \\
    &\approx& [ 0.1835, 0.2441 ].
  \end{eqnarray*}
  Note that the physical Laplacian is just a multiple of the combinatorial
  Laplacian ($\Lp = 6 \Lc$). Finally, Theorem \ref{t:KLPS} guarantees
  that there are no finitely supported eigenfunctions in $\G_{6,6}$.
\end{Example}

%%%%%%%%%%%%%%%%%%%%%%%%%%%%%%%%%%%%%%%%%%%%%%%%%%%%%%%%%%%%%%%%%%%%%%%%%%%

\section{Proof of Theorem 1}

The heart of the proof of Theorem \ref{t:cheeg} is Proposition
\ref{p:Harm} below. An earlier version of this proposition in the dual
setting (see \cite[Prop. 2.1]{BP1}) was originally obtained by helpful
discussions with Harm Derksen. Let us first introduce some
important notions related to a locally tessellating planar graph
$\G=(\V,\E,\F)$.

For a finite set $W\subseteq \V$ let $\G_W=(W,\E_W,\F_W)$ be the subgraph of
$\G$ induced by $W$, where $\E_W$ are the edges in $\E$ with both end
points in $W$ and $\F_W$ are the faces induced by the graph
$(W,\E_W)$. Euler's formula states for a finite and connected subgraph
$\G_W$ (observe that $\F_W$ contains also the unbounded face):
\begin{equation}\label{e:Euler}
  |W|-|\E_W|+|\F_W|=2.
\end{equation}
By $\dd_F W$, we denote the set of faces in $F$ which contain an edge
of $\dd_E W$. Moreover, we define the \emph{inner degree} of a face
$f\in\dd_F W$ by
$${|f|^i}_W=|f\cap W|.$$

In the following, we need the two important formulas which hold for
arbitrary finite and connected subgraphs $\G_W=(W,\E_W,\F_W)$. The first
formula is easy to see and reads as
\begin{equation}\label{e:E_W}
  \sum_{v\in W}|v|=2|\E_W|+|\dd_E W|.
\end{equation}
Since $W$ is finite, the set $\F_W$ contains at least one face which is
not in $\F$, namely the unbounded face surrounding $\G_W$, but there can
be more. Define $C(W)=|\F_W|-|\F_W\cap \F| \ge 1$. Note that $|\F_W \cap \F|$
is the number of faces in $\F$ which are entirely enclosed by edges of
$\E_W$. Sorting the following sum over vertices according to faces
gives the second formula
\begin{eqnarray}
\sum_{ v\in W}\sum_{f\ni v}\frac{1}{|f|} &=& |\F_W \cap \F| +
\sum_{f\in\dd_F W}\frac{{|f|}^i_W}{|f|} \nonumber \\
&=& |\F_W| - C(W) + \sum_{f\in\dd_F W}\frac{{|f|}^i_W}{|f|}. \label{e:F_W}
\end{eqnarray}

\begin{Prop}\label{p:Harm}
  Let $\G=(\V,\E,\F)$ be a locally tessellating planar graph and $W
  \subset \V$ be a finite set of vertices such that the induced
  subgraph $\G_W$ is connected. Then we have
  $$\ka(W)=2 -C(W) -\frac{|\dd_E W|}{2}+\sum_{f\in\dd_F W}\frac{{|f|}^i_W}{|f|}$$
\end{Prop}

\begin{proof}
  By the equations \eqref{e:E_W}, \eqref{e:F_W} and \eqref{e:Euler} we conclude
  \begin{eqnarray*}
    \ka(W)&=&\sum_{ v\in W}\ab{1-\frac{|v|}{2}+\sum_{f\ni
        v}\frac{1}{|f|}}\\
    &=&|W|-{|\E_W|}-\frac{|\dd_E W|}{2}+{|\F_W|}-C(W)+
    \sum_{f\in\dd_F
      W}\frac{{|f|}^i_W}{|f|}\\
    &=& 2 -C(W) -\frac{|\dd_E W|}{2}+\sum_{f\in\dd_F W}\frac{{|f|}^i_W}{|f|}.
  \end{eqnarray*}
\end{proof}

\begin{Prop}\label{p:dd_E W}
  Let $G=(V,E,F)$ be a locally tessellating planar graph and $3 \le q
  \le \infty$ such that $|f|\leq q$ for $f\in F$. Let $W \subset \V$
  be a finite set of vertices such that the induced subgraph $\G_W$ is
  connected. Then we have
  $$|\dd_E W|\geq\frac{2q}{q-2}(2-C(W)-\ka(W)).$$
\end{Prop}

\begin{proof}
  Since $\G$ is locally tessellating, every edge $e \in \dd_E W$
  separates precisely two different faces. The edge obtains a
  direction by its start vertex to be in $\V \backslash W$ and its end
  vertex to be in $W$. Thus it makes sense to refer to the faces at
  the left and right side of the edge $e$. Thus every edge $e \in
  \dd_E W$ determines a unique corner $(v,f) \in W \times \dd_F W$,
  where $v \in W$ is the end vertex of $e$ and $f$ is the face at the
  left side of $e$. The so defined map $\dd_E W \to W \times \dd_F W$
  is clearly injective, and thus we have
  $$
  \sum_{f\in\dd_F W} {{|f|}^i_W} =  |\{(v,f)\in W\times \dd_F W: v\in f\}|
  \ge |\dd_E W|.
  $$
  Using this fact and $|f| \leq q$ for all $f\in \F$, we conclude with
  Proposition \ref{p:Harm}
  $$ 2 - C(W) - \ka(W) = \frac{|\dd_E W|}{2} -
  \sum_{f\in\dd_F W}\frac{{|f|}^i_W}{|f|} \le |\dd_E W| \left(
    \frac{1}{2} - \frac{1}{q}\right),$$
  which proves the inequality in the proposition.
\end{proof}

Note that the Cheeger constants in Definition \ref{d:cheeg} are
obtained by taking the infimum of a particular expression over all
finite subsets $W \subset \V$. In fact, we can restrict ourselves to
consider only finite sets $W$ for which the induces graph $\G_W$ is
connected. This follows from the observation that, for a given
finite set $W \subset \V$, we can always find a non-empty subset
$W_0 \subset W$ such that $\G_{W_0}$ is a connected component of
$\G_w$ and that $|\dd_E W_0|/\vol(W_0) \le |\dd_e W|/\vol(W)$ or
$|\dd_E W_0|/|W_0| \le |\dd_e W|/|W|$, respectively. We can
reduce the sets under consideration even further. Let $W \subset \V$
be a finite set such that $\G_W$ is connected. Note that $\G_w$ has
only one unbounded face. By adding all vertices of $\V$ contained in
the union of all bounded faces of $\G_w$, we obtain a bigger finite
set $P_W \supset W$ such that $C(P_W) = 1$. (Note that all bounded
faces of $\G_{P_W}$ are also faces of the original graph $\G$.) We
call a finite set $P \subset \V$ with connected graph $\G_P$ and
$C(P) = 1$ a {\em polygon}. Clearly, we have $|\dd_E P_W |/\vol(P_W)
\le |\dd_e W|/\vol(W)$ and $|\dd_E P_W|/|P_W| \le |\dd_e W|/|W|$.
Thus it suffices for the definition of the Cheeger constants to take
the infimum only over all polygons.

With this final observation we can now prove Theorem \ref{t:cheeg}.

\begin{proof}[Proof of Theorem \ref{t:cheeg}]
Let $W \subset \V$ be a polygon. Since $C(W) = 1$, we conclude from
Proposition \ref{p:dd_E W} that
$$ \frac{|\dd_E W|}{|W|} \ge \frac{2q}{q-2} \frac{-\ka(W)}{|W|} \ge
\frac{2q}{q-2} a. $$
Taking the infimum over all polygons yields part (a) of the theorem.

For the proof of part (b), recall that $-\ka(v) \ge c \cdot |v|$ for
all vertices $v \in \V$. This implies that
$$ \frac{-\ka(W)}{\vol(W)} = \frac{- \sum_{v \in W} \ka(v)}{\sum_{v \in W}
  |v|} \ge c, $$
and, consequently, for polygons $W \subset \V$,
$$ \frac{|\dd_E W|}{\vol(W)} \ge \frac{2q}{q-2} \frac{-\ka(W)}{\vol(W)} \ge
\frac{2q}{q-2} c. $$
The statement follows now again by taking the infimum over all polygons.
\end{proof}

%%%%%%%%%%%%%%%%%%%%%%%%%%%%%%%%%%%%%%%%%%%%%%%%%%%%%%%%%%%%%%%%%%%%%%%%%%%

\section{Proof of Theorem 2}

Parts (a) and (b) of Theorem \ref{t:expgr} have very different proofs. We
present them separately.

\begin{proof}[Proof of Theorem \ref{t:expgr} (a)]
  We choose a vertex $v_0 \in \V$ and introduce the following
  functions $m, M: \F \to \{ 0,1,2,\dots, \infty \}$:
  \begin{eqnarray*}
  m(f) &=& \min \{ d(w,v_0) \mid w \in \partial f \}, \\
  M(f) &=& \max \{ d(w,v_0) \mid w \in \partial f \}.
  \end{eqnarray*}
  Note that the face $f$ ``opens up'' at distance $m(f)$ and ``closes
  up'' at distance $M(f)$ from $v_0$.  We call a face $f$ {\em finite}, if $M(f)
  < \infty$.

  The idea of the proof is to ``open up'' successively every
  finite face $f \in \F$ into an infinigon without violating the
  vertex bound. In this way, we will build up a comparison tree $\T$ with the
  same vertex bound $p$ and satisfying $\mu(\G) \le \mu(\T)$. It turns
  out, however, that finite faces $f$ with more than one vertex in the
  sphere $S_{M(f)}(v_0)$ cause problems in this ``opening up''
  procedure (since the distance relations to the vertex
  $v_0$ will be changed). Therefore, we first modify the tessellation
  $\G$ by removing all edges connecting two vertices $v,w$ at the
  same distance to $v_0$.  The modified planar graph is denoted by
  $\G_0 = (\V_0,\E_0,\F_0)$. To keep track, we add at each
  of the vertices $v,w$ a short terminal edge. These terminal edges do
  not belong ``officially'' to the graph $\G_0$ and serve merely as
  reminders that an edge can be added in their place
  without violating the vertex bound of the graph. Moreover, we can
  only guarantee $\mu(\G_0) \ge \mu(\G)$, if these inofficial edges
  are included in $\G_0$. (At the end of the procedure we will replace all
  ``inofficial'' terminal edges by infinite trees rooted in $v$ and
  $w$.) The modification $\G \to \G_0$ is illustrated in Figure
  \ref{openup1}. (For convenience, the vertices belonging to
  distance spheres $S_n(v_0)$ are arranged to lie on concentric
  Euclidean circles around $v_0$.)

  \begin{figure}[h]
    \begin{center}
      \psfrag{G}{\Large $\G$}
      \psfrag{G'}{\Large $\G_0$}
      \psfrag{v0}{$v_0$}
      \includegraphics[width=12cm]{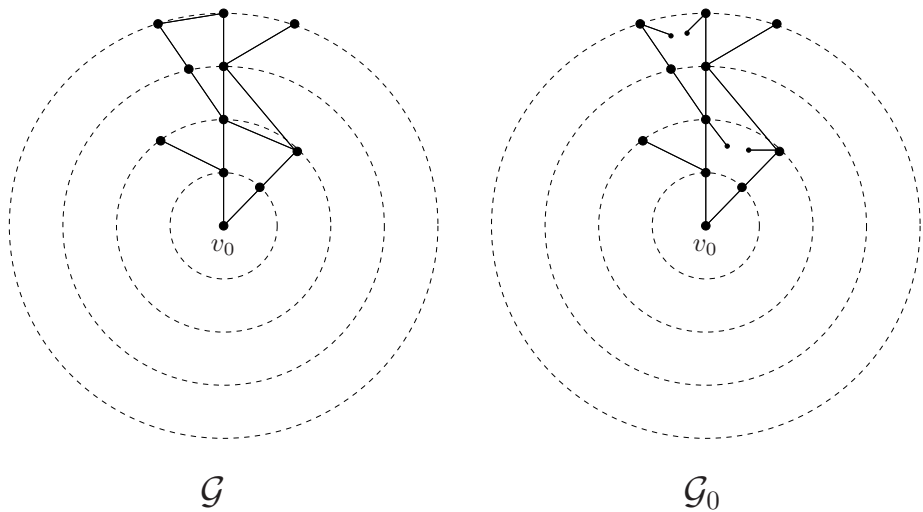}
    \end{center}
    \caption{Removing edges between vertices on the same spheres and
      replacing them by ``inofficial'' terminal edges}
    \label{openup1}
  \end{figure}

  Note that none of the distance relations of the vertices in $\G_0$
  (without the inofficial terminal edges) to the vertex $v_0$ are
  changed and that we still have $\Cut(v_0) = \emptyset$. Moreover,
  the modified graph $G_0$ (without the inofficial terminal edges) has
  a new set of faces $\F_0$. Every finite face $f$ of $\G_0$ has now
  even degree, since $f$ opens up at a single vertex in the sphere
  $S_{m(f)}(v_0)$ and $f$ closes up at a single vertex in the sphere
  $S_{M(f)}(v_0)$.

  We order all finite faces $f_0, f_1, f_2, \dots$ of $\G_0$ such
  that we have
  $$ M(f_0) \le M(f_1) \le M(f_2) \le ... $$

  Next we explain the first step of our procedure, namely, how to open
  up $f_0$ into an infinigon $\widetilde f_0$. Let $n = M(f_0) \ge 1$
  and $w \in \partial f_0$ such that $d(w,v_0) = n$. Since $C(v_0) =
  \emptyset$, we can find an infinite geodesic ray $w_0= w, w_1, w_2,
  \dots \in \V$ such that $d(w_i,v_0) = n+i$. We may think of $v_0$ as
  being the origin of the plane and of $w_0, w_1, \dots$ as being
  arranged to lie on the positive vertical coordinate axis at heights
  $n,n+1, \dots$ with straight edges between them. Now we cut our
  plane along this geodesic ray, i.e., replace the ray by two parallel
  copies of the ray and thus preventing the face $f_0$ from closing up
  at distance $n$. In this way, $f_0$ becomes an infinigon, which we
  denote by $\widetilde f_0$. (In fact, we rotationally shrink the
  angle $2\pi$ to $2\pi-\epsilon$ around $v_0$ to open up a conic
  sector of angle $\epsilon$ containing the infinigon
  $\widetilde f_0$.) The procedure is illustrated in Figure
  \ref{openup2}. Note that the vertices $w_i$ are replaced by two
  copies $w_i^{(1)}, w_i^{(2)}$, such that $w_i^{(j)}$ is connected to
  $w_{i+1}^{(j)}$ for $j=1,2$ and $w_i^{(1)}$ inherits all previous
  neighbors of $w_i$ at one side of the ray and $w_i^{(2)}$ inherits
  all previous neighbors of $w_i$ at the other side of the ray
  (this concerns in particular also the ``inofficial'' vertices). In
  this way we obtain a new planar graph $\G_1 = (\V_1,\E_1,\F_1)$.

  \begin{figure}[h]
    \begin{center}
      \psfrag{G0}{\Large $\G_0$}
      \psfrag{G1}{\Large $\G_1$}
      \psfrag{v0}{$v_0$}
      \psfrag{f0}{$f_0$}
      \psfrag{F0}{$\widetilde f_0$}
      \includegraphics[width=12cm]{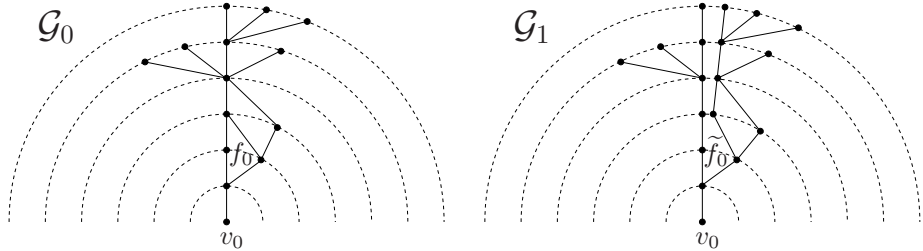}
    \end{center}
    \caption{Changing the finite face $f_0$ into an infinigon $\widetilde f_0$}
    \label{openup2}
  \end{figure}

  The graph $\G_1$ is still connected. Note also that we have
  \begin{eqnarray}
    | w_0^{(1)} | + | w_0^{(2)} | &=& | w_0 | + 1, \label{w0}\\
    | w_i^{(1)} | + | w_i^{(2)} | &=& | w_i | + 2, \quad \forall \, i \ge 1.
    \label{wi}
  \end{eqnarray}
  After including the inofficial terminal edges in the graph $\G_1$,
  we still have
  $$ | v | \le p \quad \forall\, v \in \V_1, $$
  and \eqref{w0}, \eqref{wi} imply that $\mu(\G_1) \ge \mu(\G_0) \ge \mu(\G)$.

  In the second step we carry out the same procedure with the face
  $f_1 \in \F_1$, and obtain a new connected planar graph $\G_2 =
  (\V_2,\E_2,\F_2)$, with $f_1 \in \F_1$ replaced by the infinigon
  $\widetilde f_1 \in \F_2$. Again, after including the inofficial
  terminal edges, the graph $\G_2$ has vertex bound $p$ and satisfies
  $\mu(\G_2) \ge \mu(\G_1) \ge \mu(\G)$.

  It is now clear how to repeat the procedure. Note that for every
  radius $n \ge 1$ there is a large enough $j \ge 1$ such that the
  graphs $\G_j, \G_{j+1}, \G_{j+2}, \dots$ remain unaltered inside the
  balls $B_n(\G_k,v_0)$. This fact guarantees that there is a
  well-defined limiting graph associated to the sequence $\G_j$. This
  limit is a connected tree $\T_0$ (since all faces of $\T_0$ are
  infinigons). In $\T_0$, we replace now finally the inofficial
  terminal edges by infinite trees, rooted at the corresponding proper
  vertices of the tree $\T_0$, with branching sequence
  $1,p-1,p-1,p-1,\dots$. These infinite trees can be nicely fitted
  into the infinigons to yield an infinite planar tree $\T$ with
  vertex bound $p$ and satisfying $\mu(\T) \ge \mu(\G)$. Since we
  obviously have $\mu(\T) \le \mu(\T_p) = \log (p-1)$, the
  proof of part (a) of the theorem is finished.
\end{proof}

\begin{proof}[Proof of Theorem \ref{t:expgr} (b)]
  We prove the equivalent curvature version of the statement, given in
  the remark after the theorem. Since $|f| = q < \infty$ for all faces
  $f$, $\G$ is a tessellating plane graph in the sense of \cite{BP1} and
  we have
  $$ \ka(v) = 1 - \frac{q-2}{q} |v|. $$
  Since $\{ \ka(v) \mid v \in \V \}$ is a discrete set and bounded
  from below by $-b$, we can assume, without loss of generality, that
  $-b$ is of the form $1 - \frac{q-2}{q} p$, for some integer value $p
  \ge 3$. (In fact, $p$ is the optimal upper bound on the vertex
  degree of $\G$.)

  Let $S_n, B_n$ be the combinatorial spheres and balls in $\G$ with respect
  to a reference vertex $v_0 \in \V$ and $s_n = |S_n|$. Corollary 6.4
  of \cite{BP1} states that we have
  $$ s_{n+1} - s_n = \frac{2q}{q-2} (1-\kappa(B_n)). $$
  Applying this equation twice, we derive
  $$ s_{n+2} - 2 s_{n+1} + s_n = - \frac{2q}{q-2} \kappa(S_{n+1}) \le
  \frac{2q}{q-2} b s_{n+1}. $$
  Hence we obtain the following recursion inequality
  $$ s_{n+2} \le 2 \tau s_{n+1} - s_n, \quad s_1 \le p, \ s_0 = 1, $$
  with $\tau = 1 + \frac{q}{q-2} b \ge 1$. It is easy to see that the
  sequence
  \begin{equation} \label{alphan}
  \sigma_{n+2} = 2 \tau \sigma_{n+1} - \sigma_n, \qquad \sigma_1 = p, \
  \sigma_0 = 1,
  \end{equation}
  is strictly increasing and dominates the sequence $s_n$. Moreover,
  $\sigma_n$ describes the cardinality of a sphere of radius $n$ in
  the regular tessellation $\G_{p,q}$. This implies that $\mu(\G) \le
  \mu(\G_{p,q})$.

  Now, we return to the sequence $\sigma_n$, as defined in
  \eqref{alphan}. We first consider the case $\tau > 1$. The recursion
  formula implies that
  $$ \sigma_n = u \left( \tau - \sqrt{\tau^2 -1}
  \right)^n + v \left( \tau + \sqrt{\tau^2 -1} \right)^n, $$
  with constants $u,v \in \R$ chosen in such a way that the initial
  conditions are satisfied. Since
  $$ 0 < \tau - \sqrt{\tau^2 -1} < 1, $$
  we conclude that $v \neq 0$, for otherwise we would have $\sigma_n
  \to 0$, contradicting to the fact that $\G_{p,q}$ is an infinite
  graph. Hence, $\sigma_n$ behaves asymptotically like
  $$ \sigma_n \sim v \left( \tau + \sqrt{\tau^2 -1} \right)^n, $$
  with a positive constant $v$. This, together with \eqref{muGalt} implies that
  \begin{equation} \label{mugpq}
    \mu(\G_{p,q}) = \lim_{n \to \infty} \frac{\log \sigma_n}{n} = \log
    \left( \tau + \sqrt{\tau^2 -1} \right).
   \end{equation}
  In the case $\tau = 1$, the sequence \eqref{alphan} is simply given by
  $\sigma_n = n (p-1) +1$. Linear growth of $\sigma_n$ implies that
  $\mu(\G_{p,q})=0$, which also coincides with \eqref{mugpq}.
\end{proof}

%%%%%%%%%%%%%%%%%%%%%%%%%%%%%%%%%%%%%%%%%%%%%%%%%%%%%%%%%%%%
% REFERENCES
%%%%%%%%%%%%%%%%%%%%%%%%%%%%%%%%%%%%%%%%%%%%%%%%%%%%%%%%%%%%


\begin{thebibliography}{Tho}

%\bibitem[ASW]{ASW} M. Aizenman, R. Sims, S. Warzel, \emph{Stability of
%    the absolutely continuous spectrum of random Schrödinger operators
%    on tree graphs}, Probability Theory and Related Fields 136,
%  (2006), S. 363-394

\bibitem[Al]{Al} N.\ Alon, \emph{Eigenvalues and expanders}, Theory of
  computing (Singer Island, Fla., 1984).  Combinatorica 6 (1986),
  no. 2, 83--96.

\bibitem[Ar]{Ar} N. Aronszajn, \emph{A unique continuation theorem for
    solutions of elliptic partial differential equations or
    inequalities of second order}, J. Math. Pures Appl. (9) 36 (1957),
  235--249.

\bibitem[BP1]{BP1} O.\ Baues, N.\ Peyerimhoff, {\em Curvature and
    geometry of tessellating plane graphs}, Discrete Comput.\ Geom.\
  {\bf 25} (2001), no. 1, 141-159

\bibitem[BP2]{BP2} O. Baues, N. Peyerimhoff.
\emph{Geodesics in Non-Positively Curved Plane Tessellations},
Advances of Geometry 6, no. 2, (2006), 243-263.

\bibitem[Br]{Br} R. Brooks.
\emph{A relation between growth and the spectrum of the Laplacian}, Math. Z.
178 (1981), 501-508.

\bibitem[BS]{BS} I. Benjamini, O. Schramm.
\emph{Every graph with a positive Cheeger constant contains a tree with a
positive Cheeger constant}, Geom. funct. anal. 7, (1997), 403-419.

\bibitem[Che]{Che} J. Cheeger, {\em A lower bound for the lowest
    eigenvalue of the Laplacian}, in: Problems in analysis, asymposium
  in honor of S. Bochner, Princeton Univ. Press, Princeton, (1970),
  195-199.

\bibitem[Chu]{Chu} F. Chung, \emph{Spectral graph theory}, CBMS
  Regional Conference Series in Mathematics, 92. Published for the
  Conference Board of the Mathematical Sciences, Washington, DC; by
  the American Mathematical Society, Providence, RI, 1997.

\bibitem[CGY]{CGY} F. Chung, A. Grygor'yan, S.-T. Yau, \emph{Higher
    eigenvalues and isoperimetric inequalities on Riemannian manifolds
    and graphs}, Comm. Anal. Geom.  8 (2000), no. 5, 969--1026.

%\bibitem[CFKS]{CFKS} H. L. Cycon, R. G. Froese, W. Kirsch, B. Simon.
%  \emph{Schr"odinger Operators.} Springer Verlag (1987).

%\bibitem[DeMo]{DeMo} M. DeVos, B. Mohar, \emph{An analogue of the
%      Descartes-Euler formula for infinite graphs and Higuchi's
%      conjecture}, Trans. Amer. Math. Soc.  359 (2007), no. 7,    3287--3300.

\bibitem[Do]{Do} J. Dodziuk, \emph{Difference equations, isoperimetric
    inequality and transience of certain random walks},
  Trans. Amer. Math. Soc.  284 (1984), no. 2, 787--794.

\bibitem[DKa]{DKa} J. Dodziuk, L. Karp, \emph{Spectral and function
    theory for combinatorial Laplacians}, Geometry of Random Motion,
  (R. Durrett, M.A. Pinsky ed.) AMS Contemporary Mathematics, Vol 73,
  (1988), 25-40.

\bibitem[DKe]{DKe} J. Dodziuk, W. S. Kendall, \emph{Combinatorial
    Laplacians and isoperimetric inequality}, From Local Times to
  Global Geometry, Control and Physics, (K. D. Elworthy ed.) Longman
  Scientific and Technical, (1986), 68-75.

\bibitem[DL]{DL} H. Donnelly, P. Li, \emph{Pure point spectrum and
  negative curvature for noncompact manifolds}, Duke Math. J.  46
(1979), no. 3, 497--503.

\bibitem[Fu1]{Fu1} K.\ Fujiwara, {\em Growth and the spectrum of the Laplacian
    of an infinite graph}, T{\^o}hoku Math. J. {\bf 48} (1996), 293-302

\bibitem[Fu2]{Fu2} K. Fujiwara, \emph{Laplacians on rapidly branching
    trees}, Duke Math Jour. 83, no 1, (1996), 191-202.

\bibitem[GaHuLa]{GaHuLa} S. Gallot, D. Hulin, J. Lafontaine, \emph{Riemannian
  Geometry}, Springer Verlag (1990).

\bibitem[Gri]{Gri} D. Grieser, \emph{The first eigenvalue of the
    Laplacian, isoperimetric constants, and the max flow min cut
    theorem}, Arch. Math. (Basel) 87 (2006), no. 1, 75--85.

\bibitem[Gro]{Gro} M. Gromov, \emph{Hyperbolic groups}, Essays in
  group theory, 75--263, Math. Sci. Res. Inst. Publ., 8, Springer, New
  York, 1987

\bibitem[HJL]{HJL} O. H{\"a}ggstr{\"o}m, J. Jonasson, R. Lyons,
  \emph{Explicit isoperimetric constants and phase transitions in the
    random-cluster model}, Ann. Probab. 30 (2002), no. 1, 443-473.

\bibitem[Hi]{Hi} Y. Higuchi, \emph{Combinatorial curvature for planar
    graphs}, J. Graph Theory 38 (2001), no. 4, 220--229.

\bibitem[HiShi]{HiShi} Y. Higuchi, T. Shirai, \emph{Isoperimetric constants
    of $(d,f)$-regular planar graphs}, Interdiscip. Inform. Sci. 9
  (2003), no. 2, 221-228.

%\bibitem[HLW]{HLW} S. Hoory, N. Linial, A. Wigdersin, \emph{Expander
%    graphs and their applications}, Bull. Amer. Math. Soc. (N.S.)  43
%  (2006), no. 4, 439--561.

\bibitem[Ke]{Ke} M. Keller, \emph{The essential spectrum of the
    Laplacian on rapidly branching tessellations}, arXiv:0712.3816.

\bibitem[KLPS]{KLPS} S. Klassert, D. Lenz, N. Peyerimhoff, P. Stollmann, \emph{Elliptic operators on planar graphs: unique continuation or eigenfunctions and nonpositive curvature}, Proc. Amer. Math. Soc. 134, (2005), 1549-1559.

\bibitem[KLS]{KLS} S. Klassert, D. Lenz, P. Stollmann,
  \emph{Discontinuities of the integrated density of states for random
    operators on Delone sets}, Commun. Math. Phys. 241, (2003),
  235-243.

\bibitem[KS]{KS} M. Kohmoto, B. Sutherland, \emph{Electronic states on
    a Penrose lattice}, Phys. Rev. Lett. 56, (1986), 2740-2743.

\bibitem[LV]{LV} D. Lenz, I. Veseli{\'c}, \emph{Hamiltonians on
  discrete structures: Jumps of the integrated density of states and
  uniform convergence}, arXiv:0709.2836.

\bibitem[McK]{McK} H. P. McKean, \emph{An upper bound to the spectrum
    of $\Delta$ on a manifold of negative curvature},
  J. Diff. Geom. 4, (1970), 359-366.

\bibitem[MW]{MW} B. Mohar, W. Woess, \emph{A survey on spectra of infinite graphs}, Bull. London Math. Soc. 21, (1989), 209-234.

\bibitem[St]{St} D. A. Stone, \emph{A combinatorial analogue of a
      theorem of Myers} and \emph{Correction to my paper: "A
      combinatorial analogue of a theorem of Myers"}, Illinois
    J. Math.  20 (1976), no. 1, 12--21, and Illinois J. Math.  20
    (1976), no. 3, 551--554

\bibitem[Sull]{Sull} D. Sullivan, \emph{Related aspects of positivity
  in Riemannian geometry}, J. Differential Geom.  25 (1987), no. 3,
327--351.

\bibitem[Sun]{Sun} T. Sunada, \emph{Fundamental groups and
    Laplacians}, Selected papers on number theory, algebraic geometry,
  and differential geometry, 19--32, Amer. Math. Soc. Transl. Ser. 2,
  160, Amer. Math. Soc., Providence, RI, 1994.

\bibitem[We]{We} A. Weber, \emph{Analysis of the physical Laplacian
    and the heat flow on a locally finite graph}, arXiv:0801.0812.

\bibitem[Woe]{Woe} W. Woess, \emph{A note on tilings and strong
    isoperimetric inequality}, Math. Proc. Cambridge Philos. Soc.  124
  (1998), no. 3, 385--393
\bibitem[Woe2]{Woe2} W. Woess, \emph{Random Walks on Infinite Graphs and Groups},
Cambridge Tracts in Mathematics 138, Cambridge University Press,
(2000).
\bibitem[Woj]{Woj} R. K. Wojciechowski, \emph{Stochastic Completeness
    of Graphs}, arXiv:0712.1570.

\end{thebibliography}
\end{document}